# An Efficient Analyses of the Behavior of One Dimensional Chaotic Maps using 0-1 Test and Three State Test

Joan S. Muthu,
*Department of Computer Science and Engineering,*
*SRM Institute of Science and Technology,*
*Kattankulathur, Tamil Nadu – 603203, India.*
joans@srmist.edu.in

Aditya Jyoti Paul [1,2],
*1. Department of Computer Science and Engineering,*
*SRM Institute of Science and Technology,*
*Kattankulathur, Tamil Nadu – 603203, India.*

*2. Cognitive Applications Research Lab, India.*
aditya_jyoti@srmuniv.edu.in

P. Murali
*Department of Computer Science and Engineering,*
*SRM Institute of Science and Technology,*
*Kattankulathur, Tamil Nadu – 603203, India.*
muralip@srmist.edu.in

*Abstract*— **In this paper, a rigorous analysis of the behavior of the standard logistic map, Logistic Tent system (LTS), Logistic-Sine system (LSS) and Tent-Sine system (TSS) is performed using 0-1 test and three state test (3ST). In this work, it has been proved that the strength of the chaotic behavior is not uniform. Through extensive experiment and analysis, the strong and weak chaotic regions of LTS, LSS and TSS have been identified. This would enable researchers using these maps, to have better choices of control parameters as key values, for stronger encryption. In addition, this paper serves as a precursor to stronger testing practices in cryptosystem research, as Lyapunov exponent alone has been shown to fail as a true representation of the chaotic nature of a map.**

## I. Introduction

Chaotic maps are being increasingly utilized in image cryptosystem due to their phenomenal nature of chaotic maps such as ergodicity, unpredictability, and sensitivity to their initial parameters and value [1]–[4]. The diverse disciplines where it plays significant roles are astronomy, population, security, biology, economics, social psychology and so on [5], [6] . The dynamic chaotic maps/systems are grouped into two broad categories: one dimensional (1D) and multidimensional (MD)[5]. In spite of the tediousness of the implementation and computation complexity, MD chaotic systems are preferred owing to its strong chaotic nature. On the contrary, 1D chaotic maps are simple and easy to implement. Nevertheless, it bears discontinuous range of chaotic behaviours [7]. The dynamic systems can also be classified as periodic, quasi-periodic and chaotic.

Knowledge of the strength of chaos is becoming more vital in recent times[8]. This has proportionally given rise to the development of many chaos detection techniques which are being used to analyse and measure the complex behaviour of the chaotic maps.

Several tests that are used to measure and study the chaotic behaviour are bifurcation diagram, phase portrait, largest Lyapunov exponent, Poincare map and various entropy measures with the more recent tests namely 0-1 test and three-state test (3ST) [5], [6]. Each test has its own pros and cons. Nevertheless, the largest Lyapunov exponent [9], 0-1 test [10] and 3ST [5], [6] have been proved as better analytical test.

The objective of this paper is to analyse the behaviour of 1D maps namely standard logistic map, LTS, LSS and TSS [7] using 0-1 test and 3ST. Zhou et al [7] has proved that the 1D maps exhibit a continuous chaotic behaviour in the range $r \in (0,4)$. But, frail and strong chaotic regions are detected, when these maps were analysed in a microscopic range. This result will strongly influence the selection of the range of parameter and initial values of the chaotic maps for improved results in various fields of its application.

## II. Chaos Detection Techniques

This section presents two tests namely 0-1 test and 3ST for determining the chaotic behaviour of the 1D chaotic maps.

### A. 0-1 Test

0-1 Test was proposed by Gottwald and Melbourne [11]–[12], which is used for distinguishing regular motion (periodic or quasi-periodic) from chaotic behaviour. It has been proved to be more advantages than Lyapunov exponent as the test directly applies to time series data and phase space reconstruction is not required. It gives values between 0 and 1, thereby giving a definite conclusion of the strength of chaos, whereas such a conclusion cannot be defined in largest Lyapunov exponent. Moreover, implementation of 0-1 test is more viable in non-smooth systems where linearization of the system is infeasible [9].

Definition 1 presents 0-1 test and brief details about 0-1 test is given here. For further details and better understanding of the test, the readers are advised to refer Sun et al [10].

*Definition 1*: 0-1 Test is defined by the growth rate $K$ as

$$K = \lim_{n \to \infty} \frac{\log(M(n))}{\log(n)} \quad (1)$$

where $K$ takes the value between 0 and 1, and the range of $K$ value closer to 0 depicts that the system is regular and $K$ closer to 1 depicts that the system is chaotic [9].

$M(n)$ is defined as the mean square displacement of a 2-dimensional system. It is defined as

$$M(n) = \lim_{n \to \infty} \frac{1}{N} \sum_{j=1}^{N} \left( [p(j+n) - p(j)]^2 + [q(j+n) - q(j)]^2 \right) \quad (2)$$

n=1,2,3,....

$p(n+1)$ and $q(n+1)$ is a 2-dimensional system derived from the time series $\phi(n)$ for n=1,2,3,.. and $c \in (0, 2\pi)$

$$p(n+1) = p(n) + \phi(n+1)\cos(cn),$$
$$q(n+1) = q(n) + \phi(n+1)\sin(cn) \quad (3)$$

For a time series representing regular system, the 2-dimensional system in (3) is bounded. Consequently $M(n)$ is bounded and returns $K=0$ as growth rate. The behaviour of the 2-dimensional system is similar to Brownian motion when the system is chaotic and $K$ becomes closer to 1 as the mean square displacement grows linear. But 0-1 test fails to distinguish quasi-periodic from periodic orbits [5].

### B. Three-State Test

The 3ST detects chaos based on ordinal pattern analysis and determines the period in time series [6]. The chaos indicator defined by 3ST is defined by the growth rate $K$

$$K(n) = \lim_{N \to \infty} \mu(N, n) \quad (4)$$

where $\mu(N, n)$ is the growth rate, $N$ is the length of the data series and $n$ is the smallest observation duration for the largest slope to be evaluated. The condition $n<N$ is considered for the evaluation of $K$. Here $\mu(N, n)$ is defined by the equation

$$\mu(N, n) = \frac{\log(1 + \sigma_S(N, n))}{\log N} \quad (5)$$

where $S$ is the largest slope. The measurement of the standard deviation $\sigma_S$ of S is given by

$$\sigma_S(N, n) = \sqrt{\frac{1}{Q} \sum_{j=0}^{Q-1} (S_j - \overline{S})^2} \quad (6)$$

and $\overline{S}$ is calculated with the equation

$$\overline{S} = \frac{1}{Q} \sum_{j=0}^{Q-1} S_j \quad (7)$$

where $S_j$ is the slope of the subset of the time series data. Equation (6) gives the study of the behaviour of $S$. It remains constant for periodic behaviour of the map. $S$ increases to a limiting value for quasi-periodic dynamics and it increases continuously for chaotic behavior.

### III. STANDARD AND NEW 1D CHAOTIC SYSTEMS

Logistic map is one of the most researched 1D maps owing to its simplicity in structure and implementation. This section gives an insight of the standard Logistic map and the three new 1D chaotic systems namely LTS, LSS and TSS developed by Zhou et al [7].

### A. Standard Logistic Map

The standard logistic map is a simple dynamical system which exhibits phenomenal chaotic behavior. It is defined by the mathematical definition

$$X_{n+1} = rX_n(1 - X_n) \quad (8)$$

where $r \in (0,4]$. The logistic map exhibits chaotic behavior in the range of $r \in [3.57, 4]$. Hence, the range of the chaotic behavior is very limited. Moreover, there is no continuous chaotic behavior in this range which is evident from the periodic windows in the bifurcation diagram as shown in Figure 1 (a). This has limited the application of logistic map and consequently given rise to devising of new modified logistic maps with improved chaotic behavior.

### B. Logistic-Tent system

LTS was designed from the seed maps Logistic map and Tent map [7]. It is proved that this map displays better chaotic behavior with wider range i.e. for $r \in (0,4]$. The LTS is defined by the mathematical expression

$$X_{n+1} = \begin{cases} (rX_n(1-X_n) + (4-r)X_n/2) \bmod 1 & X_i < 0.5 \\ (rX_n(1-X_n) + (4-r)(1-X_n)/2) \bmod 1 & X_i \geq 0.5 \end{cases} \quad (9)$$

where $r \in (0,4]$. This map exhibits uniform distribution of output sequence in bifurcation diagram as shown in Figure 1 (b).

### C. Logistic-Sine system

The LSS is defined by the mathematical expression given in Equation (10).

$$X_{n+1} = (rX_n(1-X_n) + (4-r)\sin(\pi X_n)/4) \bmod 1 \quad (10)$$

where $r \in (0,4]$. The bifurcation diagram in Figure 1 (c) shows a wider chaotic behaviour in the range $r \in (0,4]$.

## D. Tent-Sine system

The mathematical expression of TSS [7] is given by equation (11).

$$X_{n+1} = \begin{cases} (rX_n/2 + (4-r)\sin(\pi X_n)/4) \bmod 1 & X_i < 0.5 \\ (r(1-X_n)/2 + (4-r)\sin(\pi X_n)/4) \bmod 1 & X_i \geq 0.5 \end{cases}$$
(11)

where $r \in (0,4]$. The bifurcation diagram in Figure 1 (d) depicts a uniform range.

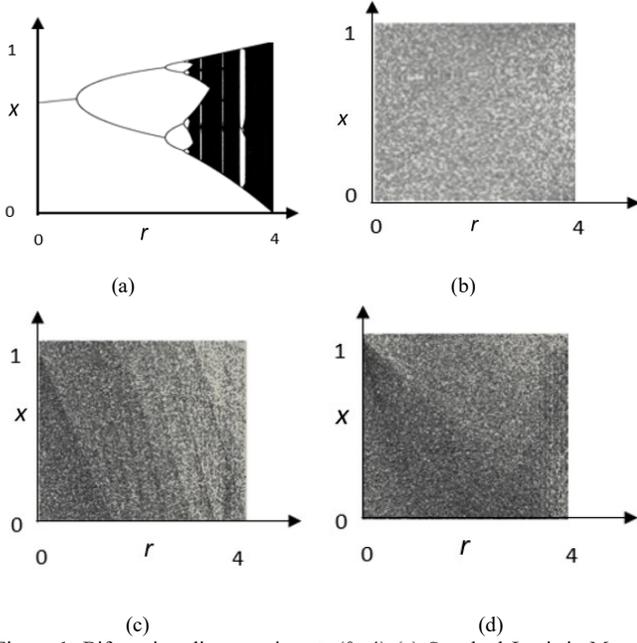

Figure 1. Bifurcation diagrams in $r \in (0, 4)$ (a) Standard Logistic Map (b) Logistic-Tent system (c) Logistic-Sine system and (d) Tent-Sine system

## IV. RESULTS AND DISCUSSION

0-1 Test and 3ST was performed on the standard logistic map, LTS, LSS and TSS in the range $r \in (3.1, 4)$.

### A. Analyses using 0-1 Test

0-1 test was experimented on all the above mentioned 1D maps with parameters $N=5000$, $x_0=0.01$ and $c=0.8$. For this, the time series $\phi(n)$ is obtained for the logistic map, LTS, LSS and TSS by iterating the equations (8)-(11) individually. Then, the 2-dimensional system of $p(n+1)$ and $q(n+1)$ of Equation (3) is evaluated from time series $\phi(n)$ for each map. Figure (2) presents the $p$-$q$ graph of the logistic map for $r=3.15$ and $r=3.95$ respectively.

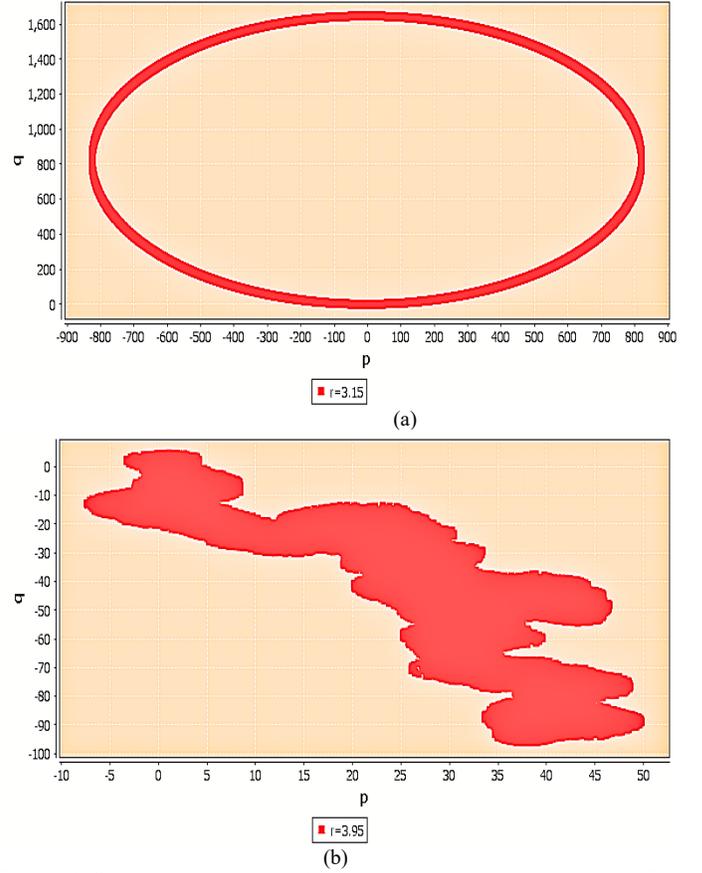

Figure 2. Standard Logistic map (a) 2-dimensional system is bounded for $r=3.15$ (b) 2-dimensional system displays Brownian motion for $r=3.95$

Figure 2 (a) depicts a bounded 2-dimensional system represented by the Equation (3) for $r=3.15$ of the logistic map, which represents the presence of a regular dynamics. While, at $r=3.95$ of the logistic map, the $p$-$q$ graph exhibits Brownian motion in Figure 2 (b), symbolizing the chaotic behavior.

Figure 3 indicates the corresponding slope of growth rate of $K$ of the logistic map at $r=3.15$ and $r=3.95$.

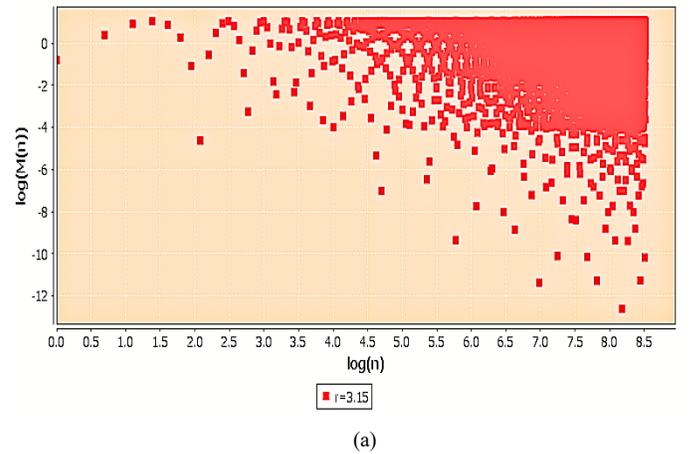

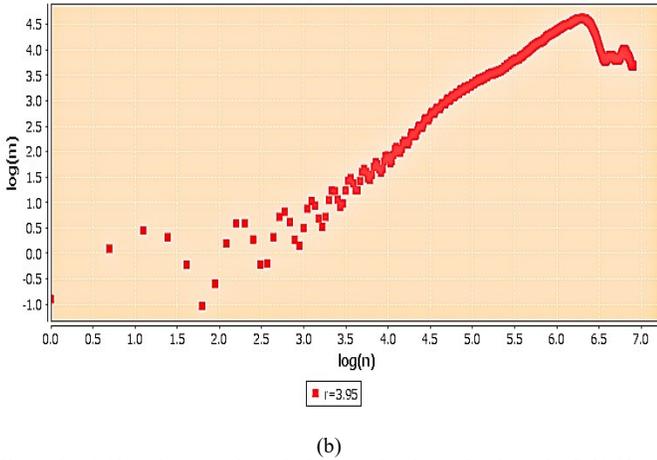

(b)

Figure 3. (a) Slope is approximately 0 in regular dynamics for *r*=3.15 (b) Slope is approximately 1 for chaotic dynamics for *r*=3.95

Figure 3 (a) indicates the slope of growth rate of *K* is approximately 0 for *r*=3.15 and Figure 3 (b) is indicative of the occurrence of chaos at *r*=3.95 as the slope approaches 1. Hence, the 0-1 test manifests the regular and chaotic dynamics in the logistic map.

The dynamic nature of the 1D chaotic maps LTS, LSS and TSS is also experimented and the 0-1 test is performed on these maps for $r \in (3.1, 4)$. The *K* values obtained for the listed *r* values are given in Table I.

TABLE I. 0-1 TEST RESULTS: *K* VALUES OBTAINED FOR THE R VALUES MENTIONED

| 1D chaotic maps | Standard Logistic map | Logistic Tent system | Logistic Sine System | Tent Sine System |
|---|---|---|---|---|
| r | K | K | K | K |
| 3.15 | 0.0482 | 0.6054 | 0.8335 | 0.7709 |
| 3.25 | 0.0468 | 0.5686 | 0.7394 | 0.6857 |
| 3.35 | 0.0455 | 0.6758 | 0.7692 | 0.7700 |
| 3.45 | 0.0444 | 0.6992 | 0.6158 | 0.7641 |
| 3.55 | 0.0579 | 0.6570 | 0.6739 | 0.6516 |
| 3.65 | 0.1921 | 0.7513 | 0.7112 | 0.6306 |
| 3.75 | 0.4119 | 0.7004 | 0.6438 | 0.7121 |
| 3.85 | 0.0094 | 0.5172 | 0.5503 | 0.4176 |
| 3.95 | 0.6098 | 0.4678 | 0.5137 | 0.0847 |

Table I concludes that though the 0-1 test shows slope approaching 1 for all values in the range $r \in (3.1, 4)$ for LTS, LSS and TSS, there is a variation in the *K* value, depicting weak and strong chaotic regions. Hence, these maps which are claimed to be chaotic in the range (0,4) [7] does not possess uniform chaotic nature throughout the specified range.

It further projects that LSS possesses the strongest chaotic nature. Nevertheless, it is emphasized that the *K* value of the 0-1 test could not differentiate periodic and quasi-periodic behaviour of a map.

*B. Analyses using 3ST*

We observed that 3ST is a more accurate test to determine the chaotic nature of a map as it distinguished periodic, quasi-periodic and chaotic behaviour [6]. Using the algorithm described in the work of Eyebe and Koepf [6], the 3ST was tested with the following parameters: *N*=5000, *n*=50, *Q*=100, *p*=1. The 3ST is performed on the logistic map, LTS, LSS and TSS in the range $r \in (3.1, 4)$. Surprisingly, it revealed three behaviour at different *r* values, viz. periodic, quasi-periodic and chaotic nature.

Figure 4 shows the results obtained for $r \in (2, 2.9)$ and $r \in (3.6, 3.69)$ of LTS.

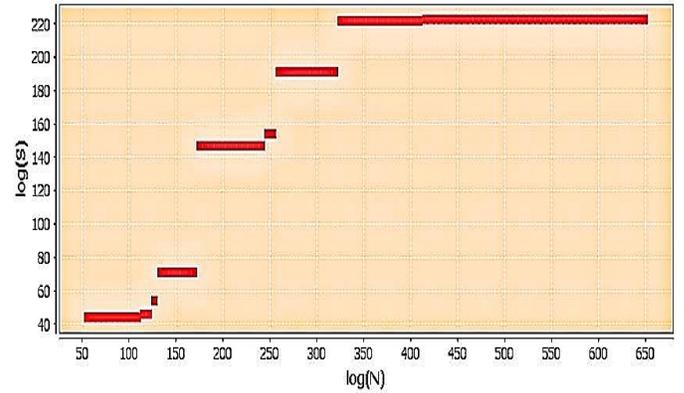

(a)

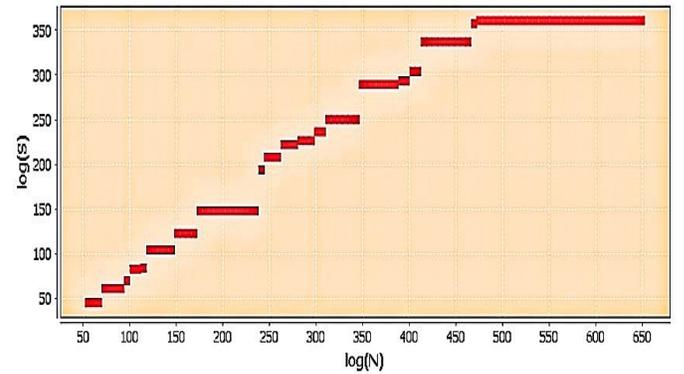

(b)

Figure 4. 3ST (a) S increases to a limiting value for quasi-periodic dynamics for $r \in (2, 2.9)$ of LTS (b) S increases continuously for chaotic behavior $r \in (3.6, 3.69)$ of LTS

It is observed that Figure 4 (a) shows large steps indicating the quasi periodic behaviour in $r \in (2, 2.9)$ of LTS. While, the steep slope in Figure 4 (b) represents the chaotic behaviour in $r \in (3.6, 3.69)$. Similarly, the results of 3ST for the other ranges of *r* was evaluated for all the maps. The comparison of the results of 3ST is tabulated in Table II.

TABLE II. COMPARISON OF THE 3ST RESULTS INDICATING THE BEHAVIOR OF THE 1D MAPS

| 1D chaotic maps $r$ | | Standard Logistic Map | Logistic Tent System | Logistic Sine System | Tent Sine System |
|---|---|---|---|---|---|
| 3.1 | 3.19 | Periodic | quasi-periodic | chaotic | chaotic |
| 3.2 | 3.29 | Periodic | quasi-periodic | chaotic | quasi-periodic |
| 3.3 | 3.39 | Periodic | chaotic | chaotic | chaotic |
| 3.4 | 3.49 | Periodic | chaotic | quasi-periodic | chaotic |
| 3.5 | 3.59 | quasi-periodic | chaotic | quasi-periodic | quasi-periodic |
| 3.6 | 3.69 | quasi-periodic | chaotic | chaotic | quasi-periodic |
| 3.7 | 3.79 | Chaotic | chaotic | quasi-periodic | chaotic |
| 3.8 | 3.89 | Periodic | quasi-periodic | quasi-periodic | quasi-periodic |
| 3.9 | 3.99 | Chaotic | quasi-periodic | quasi-periodic | quasi-periodic |

TABLE III. COMPARISON OF THE BEHAVIOUR OF THE 1D CHAOTIC SYSTEMS

| 1D Map | Chaotic | | Quasi-periodic | | Periodic | |
|---|---|---|---|---|---|---|
| Logistic map | 3.7 | 3.79 | 3.5 | 3.69 | 3.1 | 3.49 |
| | 3.9 | 3.99 | - | - | 3.8 | 3.89 |
| Logistic Tent system | 3.3 | 3.79 | 3.1 | 3.29 | - | - |
| | - | - | 3.8 | 3.99 | - | - |
| Logistic Sine System | 3.1 | 3.39 | 3.4 | 3.59 | - | - |
| | 3.6 | 3.69 | 3.7 | 3.99 | - | - |
| Tent Sine System | 3.1 | 3.19 | 3.2 | 3.29 | - | - |
| | 3.3 | 3.49 | 3.5 | 3.69 | - | - |
| | 3.7 | 3.79 | 3.8 | 3.99 | - | - |

Needless to say, Table II clearly categories three dynamic behaviours and distinguishes quasi-periodic motions from periodic and chaotic motions. Further, it is observed that the chaotic behaviour is not evenly distributed and some regions are found to be quasi-periodic in LTS, LSS and TSS. Hence, the claim in [7] that LTS, LSS and TSS are chaotic throughout is debatable. It is also noted from Table II that LTS has wider chaotic region than other maps.

*C. Discussion*

The new 1D chaotic systems were proved to have positive Lyapunov exponent for the range $r \in (0,4]$ in the work of Zhou et al [7]. Further, output sequence maps the entire range between 0 and 1.

In this work, the behaviour of the logistic map, LTS, LSS and TSS are examined in $r \in (3.1, 4)$ with the 0-1 test and 3ST. The 0-1 test proved different strengths of chaos in these maps. The 3ST is another remarkable test that classified the dynamic motions as periodic, quasi-periodic and chaotic. These behaviours are summarized in Table III, which depicts the dynamics in the various region of $r$.

Table III clearly portrays stronger and weaker regions of chaos and it is clearly evident that strong chaotic regions is limited in the region specified, as some regions exhibit quasi-periodic behaviour. Nonetheless to mention, LTS, LSS and TSS do not exhibit periodic motions, but exhibit quasi-periodic motion in some regions of test.

V. CONCLUSION

In this work, a microscopic analysis is performed on the proposed maps with 0-1 test and 3ST to examine the strength of chaos exhibited in these maps at $r \in (3.1, 4)$. 0-1 test demonstrated that LTS, LSS and TSS are not regular at all the values of $r \in (3.1, 4)$. But, it exhibited different strengths of chaos throughout the $r$ value. The experimental results showed that LSS possesses the strongest chaotic nature.

Next, the 3ST was performed on the proposed maps. Interestingly, the analysis showed quasi-periodic behaviour at intermediate regions in all of these maps. The 3ST results proved that LTS has wider chaotic range than LSS and TSS.

Hence, it is concluded that LTS, LSS and TSS did not exhibit a uniform chaotic nature and the results show strong contradiction to the statement of Zhou et al. Further, it is to be noted that though Lyapunov exponent has the advantage of simplicity in implementation in smooth 1D chaotic systems, it fails to detect strong and weak chaotic regions. Hence, it is not a sufficient test to detect chaos in a dynamic system. 0-1 test and 3ST throws light on the strength of chaotic maps that varies across different values of $r$.

In this work, regions of weak chaos are exhibited that may not be suitable for applications requiring strong chaotic values. Further, we also present chaotic region that we recommend researchers to apply in their applications for better results. Hence, we strongly convey researchers that the selection of the range of control parameter $r$ value must be tested using 0-1 test and 3ST also for better results in different applications.